\documentclass[11pt]{amsart}
\usepackage{amsmath}
\usepackage{amsxtra}
\usepackage{amscd}
\usepackage{amsthm}
\usepackage{amsfonts}
\usepackage{amssymb}
\usepackage{eucal}

\begin{document}
\newcommand{\ga}{{\mathfrak a}}
\newcommand{\gb}{{\mathfrak b}}
\newcommand{\gn}{{\mathfrak n}}
\newcommand{\go}{{\mathfrak o}}
\newcommand{\g}{{\mathfrak{g}}}
\newcommand{\slt}{\mathfrak{sl}_2}
\newcommand{\slth}{\widehat{\mathfrak{sl}}_2}
\newcommand{\sltr}{\mathfrak{sl}_3}
\newcommand{\nub}{{\boldsymbol \nu}}
\newcommand{\mub}{{\boldsymbol \mu}}
\newcommand{\sub}{{\boldsymbol \sigma}}
\newcommand{\nb}{\mathbf{n}}
\newcommand{\mb}{\mathbf{m}}
\newcommand{\bm}{\bar{m}}
\newcommand{\bmb}{\overline{\mathbf{m}}}
\newcommand{\bu}{\bar{u}}
\newcommand{\ub}{\mathbf{u}}
\newcommand{\bub}{\bar{\mathbf{u}}}
\newcommand{\vb}{\mathbf{v}}
\newcommand{\wb}{\mathbf{w}}
\newcommand{\bv}{\bar{v}}
\newcommand{\bvb}{\bar{\mathbf{v}}}
\newcommand{\Mb}{\mathbf{M}}
\newcommand{\bMb}{\overline{\mathbf{M}}}
\newcommand{\bM}{\overline{M}}
\newcommand{\blam}{\bar{\lambda}}
\newcommand{\bmu}{\bar{\mu}}
\newcommand{\bone}{\mathbf{1}}

\newcommand{\et}{\tilde{e}}
\newcommand{\htil}{\tilde{h}}
\newcommand{\su}{\sigma}
\newcommand{\vpi}{\varpi}
\newcommand{\bpi}{\overline{\pi}}
\newcommand{\bPi}{\overline{\Pi}}
\newcommand{\bare}{\bar{e}}
\newcommand{\bga}{\bar{\ga}}

\newcommand{\bk}{\bar{k}}
\newcommand{\bp}{\bar{p}}

\newcommand{\nn}{\nonumber}
\newcommand{\bea}{\begin{eqnarray}}
\newcommand{\ena}{\end{eqnarray}}
\newcommand{\be}{\begin{eqnarray*}}
\newcommand{\en}{\end{eqnarray*}}

\newcommand{\res}{{\mathop{\rm res}}}
\newcommand{\id}{{\rm id}}
\newcommand{\ch}{{\rm ch}}
\newcommand{\End}{\mathop{{\rm End}}}
\newcommand{\tr}{\mathop{{\rm tr}}}
\newcommand{\gr}{\mathop{{\rm gr}}}
\newcommand{\ggr}{{\rm gr}}
\newcommand{\bra}[1]{\langle #1 |}        
\newcommand{\ket}[1]{{| #1 \rangle}}      
\newcommand{\br}[1]{{\langle #1 \rangle}}  
\newcommand{\qbin}[2]{{\left[
\begin{matrix}{\displaystyle #1}\\
{\displaystyle #2}\end{matrix}
\right]
}}
\newcommand{\bep}{\bar{\epsilon}}
\newcommand{\vep}{\varepsilon}
\newcommand{\ep}{\epsilon}
\newcommand{\pit}{\Pi}
\newcommand{\ut}{v}
\newcommand{\xt}{\tilde{x}}
\newcommand{\zz}{{\mathcal Z}}
\newcommand{\z}{\zeta}
\newcommand{\cK}{{\mathcal K}}
\newcommand{\cH}{{\mathcal H}}
\newcommand{\cG}{{\mathcal G}}
\newcommand{\fusn}{\circledast}
\newcommand{\wg}{\widehat{\g}}
\newcommand{\ver}{{\mathcal V}^{(k)}}
\newcommand{\usl}{U(\slt[t])}
\newcommand{\gY}{{\mathfrak Y}}
\newcommand{\gYb}{\overline{\mathfrak Y}}
\newcommand{\bY}{\overline{Y}}
\newcommand{\bX}{\overline{X}}
\newcommand{\bW}{\overline{W}}
\newcommand{\C}{{\mathbb C}}
\newcommand{\Z}{{\mathbb Z}} 
\newcommand{\N}{{\mathbb N}}
\newcommand{\R}{{\mathbb R}}
\newcommand{\Q}{{\mathbb Q}} 
\newcommand{\F}{\mathcal F}
\newcommand{\Ft}{\mathcal G}
\numberwithin{equation}{section}
\newtheorem{thm}{Theorem}[section]
\newtheorem{prop}[thm]{Proposition}
\newtheorem{lem}[thm]{Lemma}
\newtheorem{cor}[thm]{Corollary}
\newtheorem{rem}[thm]{Remark}
\newtheorem{dfn}{theorem}
\renewcommand{\theequation}{\thesection.\arabic{equation}}
\def\theenumi{\roman{enumi}}
\def\labelenumi{(\theenumi)}

\title[Two character formulas
for $\slth$ spaces of coinvariants]
{Two character formulas for $\slth$ spaces of coinvariants}
\author{B. Feigin, M. Jimbo, S. Loktev and T. Miwa}
\address{BF: Landau institute for Theoretical Physics, Chernogolovka,
142432, Russia}\email{feigin@feigin.mccme.ru}  
\address{MJ: Graduate School of Mathematical Sciences, The 
University of Tokyo, Tokyo 153-8914, Japan}\email{jimbomic@ms.u-tokyo.ac.jp}
\address{SL:  Institute for Theoretical and Experimental Physics,
B. Cheremushkinskaja, 25, Moscow 117259, Russia}
\address{Independent University of Moscow, B. Vlasievsky per, 11, 
Moscow 121002, Russia}\email{loktev@mccme.ru}
\address{TM: Division of Mathematics, Graduate School of Science, 
Kyoto University, Kyoto 606-8502
Japan}\email{tetsuji@kusm.kyoto-u.ac.jp}

\date{\today}

\begin{abstract}
We consider  $\slth$ spaces of coinvariants
with respect to two kinds of ideals of 
the enveloping algebra $U(\slt\otimes\C[t])$.
The first one is generated by $\slt\otimes t^N$, 
and the second one is generated by 
$e\otimes P(t), f\otimes \overline{P}(t)$ 
where $P(t),\overline{P}(t)$ are fixed generic polynomials.  
(We also treat a generalization of the latter.)
Using a method developed in our previous paper, 
we give new fermionic formulas for their Hilbert polynomials 
in terms of the level-restricted Kostka polynomials and 
$q$-multinomial symbols. 
As a byproduct, we obtain a fermionic formula 
for the fusion product of $\sltr$-modules with 
rectangular highest weights, generalizing a known result for 
symmetric (or anti-symmetric) tensors. 
\end{abstract}
\maketitle

\setcounter{section}{0}
\setcounter{equation}{0}
\section{Introduction}\label{sec:1}
Spaces of conformal coinvariants 
are the central objects in conformal field theory. 
In addition to their many intriguing features,   
they present also 
some combinatorial problems which are 
worth being pursued in their own right \cite{FL}. 
A systematic study of this aspect has been launched in a series of papers
\cite{FKLMM1}--\cite{FKLMM3}. 

To be specific, the spaces of coinvariants
considered in these works are the quotient spaces 
$L^{(k)}_l/\ga L^{(k)}_l$, where 
$L^{(k)}_l$ is a level $k$ integrable $\slth$-module 
and $\ga$ is a subalgebra of $\slt\otimes\C[t]$ in one of
the following forms:
\bea
&&
\ga_N=\slt\otimes Q(t)\C[t],
\label{BC}\\
&&
\ga^{(p,\bar{p})}=\C e \otimes P(t) \C[t]
\oplus \C h \otimes P(t)\overline{P}(t) \C[t]
\oplus \C f \otimes \overline{P}(t) \C[t]. 
\label{MC}
\ena
Here $e,h,f$ are the standard generators of $\slt$, and 
$Q(t)=\prod_{i=1}^N(t-\z_i)$, 
$P(t)=\prod_{j=1}^p(t-\z_j)$, 
$\overline{P}(t)=\prod_{j=1}^{\bar{p}}(t-\z_{j+p})$ 
are polynomials. 
(We use the convention that $\slt\otimes t^{-1}\C[t^{-1}]$
annihilates the highest weight vector. 
See subsection \ref{subsec:2.1}.)
These spaces of coinvariants 
are finite dimensional vector spaces equipped with  
a natural filtration by the degree in $t$. 
One of the basic quantities of interest is 
the Hilbert polynomial, or {\it character}, 
of their associated graded spaces. 
Often the characters are 
written in certain specific form called {\it fermionic formulas}. 
Their actual form depends on the method used to obtain them,  
and the same quantity can have different expressions.
Indeed, in \cite{FJKLM2}, we found  a 
fermionic formula for the case \eqref{MC} with $\bar{p}=0$
which is 
different from the one obtained earlier in \cite{FKLMM3}. 
The purpose of this paper is to present analogous formulas
in the cases \eqref{BC},\eqref{MC} and their generalizations.

The formula mentioned above has its origin 
in the dimension formula 
in terms of the Verlinde algebra $\ver$. 
Recall that $\ver$ is an associative unital ring over $\Z$, 
with basis $[l]$ ($0\le l\le k$) and multiplication rule 
\bea
[l]\cdot [l']=\sum_{|l-l'|\le l''\le \min(l+l',2k-l-l') 
\atop l''\equiv |l-l'|\bmod 2}
[l''].
\label{verlinde}
\ena
{}For an element $a=\sum_{l=0}^ka_l[l]\in\ver$ we write 
the coefficient $a_l$ as $(a:[l])_k$. 
The general dimension formula \cite{FJKLM1} applied to 
\eqref{MC} with $\bar{p}=0$ 
tells that, when $\z_i$ are distinct, we have
\be
\dim L_l^{(k)}/\ga^{(p,0)}L_l^{(k)}
&=&
\left(\bigl([0]+[1]+\cdots+[k]\bigr)^p:[l]\right)_k
\\
&=&
\sum_{m_0,\cdots,m_k\ge 0\atop m_0+\cdots+m_k=p}
\left(p \atop m_0,\cdots,m_k\right) 
([0]^{m_0}[1]^{m_1}\cdots[k]^{m_k}:[l])_k,
\en
where the second line is simply a result of expansion in the first. 
The fermionic formula obtained in \cite{FJKLM2} 
is a $q$-analog of the right hand side, 
in which the multinomial coefficient is replaced by 
the $q$-multinomial coefficient, 
and $([0]^{m_0}[1]^{m_1}\cdots[k]^{m_k}:[l])_k$ 
by its $q$-analog called the restricted Kostka polynomial. 
In fact, \cite{FJKLM2} deals with a slightly 
more general situation, wherein $([0]+\cdots+[k])^p$ is replaced
by $\prod_{j=1}^p([0]+\cdots+[k_j])$ for arbitrary
$k_1,\cdots,k_p\in\{0,\cdots,k\}$. 
The corresponding coinvariants arise by replacing Lie subalgebras by 
appropriate right ideals of the enveloping algebra \cite{FJKLM1}. 
In this paper we also treat such a generalization for \eqref{MC}. 
It is proved in \cite{FKLMM1},\cite{FKLMM2} that the characters
of the coinvariants for \eqref{BC}, \eqref{MC} do not depend on the $\z_i$.
In view of this result, 
we will restrict to the case $Q(t)=t^N$ for \eqref{BC}. 
As for the generalization of the case \eqref{MC} mentioned above, 
the independence of the character has not yet been established. 
{}For that reason, we assume for \eqref{MC} that 
$\z_1,\cdots,\z_{p+\bar{p}}$ are pairwise distinct. 

In order to obtain the fermionic formula,  
we follow the method used in \cite{FJKLM2}. 
It is rather indirect and  
consists of several steps. 
{}First we apply the general equivalence theorem in \cite{FJKLM1}
and reduce the problem to that of a fusion product 
of finite dimensional $\slt$-modules and their quotient spaces. 
Each constituent of the fusion product is 
a reducible $\slt$-module whose cyclic vector is given as a sum.
Technically it is difficult to handle such fusion products directly. 
To circumvent this point, 
we embed $\slt$ into some larger Lie algebra $\g$  
and show that the fusion product in question is the same as that of 
irreducible $\g$-modules with lowest vector as cyclic vector. 
{}For \eqref{MC} 
$\g=\sltr$, and for \eqref{BC} $\g=\go_5$. 
The final step is a reduction 
to the fusion product of modules over an abelian Lie algebra. 

Though the applicability of this method is limited,  
the resulting formulas seem to have a universal nature. 
It is left as an interesting open problem  
to establish such formulas in full generality. 

The plan of the paper is as follows.
In section \ref{sec:2} we formulate the problem and 
state the main results. 
In section \ref{sec:3} we give a proof for the case of the 
coinvariants with respect to \eqref{BC}. 
The proof in the other case \eqref{MC} goes quite parallel, 
and is briefly described in section \ref{sec:4}. 
We also obtain a fermionic formula for the 
fusion product of $\sltr$-modules involving both 
symmetric and anti-symmetric tensors (Theorem \ref{thm:sl3}).

\section{Statement of the result}\label{sec:2}
\subsection{Notation}\label{subsec:2.1}
{}First we set up the notation. 
We use the symbol
$e_{ab}=\left(\delta_{ai}\delta_{jb}\right)$ 
for matrix units, 
whose size should be clear from the context. 
Let $e,f,h$ be the standard generators of $\slt$. 
Let $L_l^{(k)}$ be the 
integrable module over the affine Lie algebra
$\slth=\slt\otimes\C[t,t^{-1}]\oplus\C K\oplus\C D$, 
with level $k$ and 
highest weight $l$ ($l,k\in \Z$, $0\le l\le k$). 
We use the following convention:
the canonical central extension is given by 
\be
[x_m,y_n]=[x,y]_{m+n}-m\delta_{m+n,0}(x|y)K
\qquad (x,y\in\slt),
\en
where
$x_m=x\otimes t^m$, $(x|y)=\mathop{\rm tr}(xy)$, and   
$K$ acts as $k$ times the identity. 
The highest weight vector 
$v_l^{(k)}\in L_l^{(k)}$ satisfies 
$x_iv_l^{(k)}=0$ ($i<0$, $x=e,f,h$), $h_0v_l^{(k)}=lv_l^{(k)}$ and
$e_1^{k-l+1}v_l^{(k)}=0$. 
In the terminology of \cite{FJKLM2}, 
$L_l^{(k)}=L_l^{(k)}(\infty)$ is `placed at infinity'. 
In general, for a Lie algebra $\g$, we write
$\g[t]=\g\otimes \C[t]$. 

Let
$V=\oplus_{\alpha\in I}V_{\alpha}$
be a graded vector space indexed by 
$I=\Z^{m+1}_{\ge 0}+\beta$, $\beta\in\Q^{m+1}$.
Let $D_i$ ($0\le i\le m$) be 
the degree operator which acts on 
$V_{\alpha}$ as $\alpha_i$ times the identity. 
We call
\be
\ch_{z_0,\cdots,z_m}V
&=&
{\rm tr}_V\Bigl(z_0^{D_0}\cdots z_m^{D_m}\Bigr)
\\
&=&\sum_\alpha z^{\alpha} \dim V_{\alpha}
\qquad
(z^\alpha=z_0^{\alpha_0}\cdots z_m^{\alpha_m})
\en
its character. 
We use the inequality sign of formal series 
$\sum_{\alpha}f_{\alpha}z^\alpha\le
\sum_{\alpha}g_{\alpha}z^\alpha$ 
to mean that $f_{\alpha}\le g_{\alpha}$ for all $\alpha$. 
We regard 
$L_l^{(k)}$ as a graded $\slt[t]$-module 
by the degree operators $D_0=td/dt$ and $D_1=h_0/2$. 
This is equivalent to the assignment 
\bea
&&\deg e_i=(i,1), 
\quad \deg f_i=(i,-1),
\quad \deg h_i=(i,0),
\label{degree}\\
&&\deg v_l^{(k)}=(0,l/2).
\nn
\ena

\subsection{Spaces of coinvariants} \label{subsec:2.2}
{}For a module $V$ over an algebra $A$ and a right ideal 
$Y\subset A$, 
we call $V/YV$ the space of coinvariants and 
use the abbreviated notation $V/Y$. 
In this paper we consider spaces of coinvariants of $L_l^{(k)}$
with respect to two kinds of ideals of $U(\slt[t])$.

The first ideal is generated by polynomial currents 
which have $N$-fold zeroes at the origin:
\be
B_N=\ga_N U(\slt[t]),\qquad 
\ga_N=\slt\otimes t^N\C[t].   
\en
The second ideal is defined as follows.
{}For $m\in\Z_{\ge 0}$, define right ideals of $U=U(\slt[t])$ by 
\be
Y_m=e_0^{m+1}U+f_0U+B_1,
\qquad
\bY_m=f_0^{m+1}U+e_0U+B_1.
\en
{}Fix distinct complex numbers 
$\zz=(\z_1,\cdots,\z_{p+\bp})$.
Choose integers
$k_1,\cdots,k_p$, 
$\bk_1,\cdots,\bk_{\bp}\in \{0,1,\cdots,k\}$,  
and set 
\bea
&\Mb=(M_1,\cdots,M_k),\quad
&M_a=\sharp\{i\mid k_i=a\},
\label{Mb1}\\
&
\bMb=(\bM_1,\cdots,\bM_k),
\quad 
&\bM_a=\sharp\{i\mid \bk_i=a\}. 
\label{Mb2}
\ena
With the above data, we associate 
the fusion right ideal in the sense of \cite{FJKLM1}, 
\bea
Y_{\Mb,\bMb}(\zz)=
Y_{k_1}\fusn\cdots\fusn Y_{k_p}
\fusn \bY_{\bk_1}\fusn\cdots\fusn \bY_{\bk_{\bp}}(\zz). 
\label{Ym}
\ena

Notation being as above, 
we consider the spaces of coinvariants
\bea
&&L_l^{(k)}/B_N,
\label{bigc}
\\
&&L_l^{(k)}/Y_{\Mb,\bMb}(\zz).
\label{mixc}
\ena
The first space \eqref{bigc} was introduced in \cite{FL} 
for general non-twisted affine Lie algebras, and was 
studied in \cite{FKLMM1} for $\slth$.   
The second space \eqref{mixc} 
appears as Example 4 in \cite{FJKLM1}. 
In the special case $M_{i}=\delta_{ik}M$ and 
$\bM_{i}=\delta_{ik}\bM$, \eqref{mixc}  
coincides with the space of coinvariants 
studied in \cite{FKLMM2,FKLMM3}.

Since $B_N$ is a homogeneous ideal, 
the space \eqref{bigc} inherits a natural bi-grading
from $L_l^{(k)}$. 
On the other hand, 
the ideal $Y_{\Mb,\bMb}(\zz)$ is not homogeneous,
and the space \eqref{mixc} is only filtered. 
Instead we consider the associated graded space
$\gr\Bigl(L_l^{(k)}/Y_{\Mb,\bMb}(\zz)\Bigr)$.
The aim of this paper is to find an expression for 
the characters of these graded spaces.

\subsection{$q$-multinomials 
and restricted Kostka polynomials}\label{subsec:2.3}
The character formulas we are going to present 
consists of three pieces:
certain $q$-multinomial symbols, restricted Kostka polynomials, 
and the character of the fusion product of $\slt$-modules. 
In this subsection we recall them.

Let 
\be
&&\qbin{m}{n}=\begin{cases}
\displaystyle\frac{(q^{m-n+1})_n}{(q)_n} & (0\le n\le m),\\
0 & \mbox{otherwise},\\
\end{cases}
\en
denote the $q$-binomial symbol, 
wherein $(z)_n=\prod_{j=0}^{n-1}(1-zq^j)$. 
{}For an array of non-negative integers 
$\mb=(m_1,\cdots,m_k)\in\Z_{\ge 0}^k$, 
we set 
\be
&&|\mb|=\sum_{i=1}^k i m_i, 
\\
&&\lambda_a(\mb)=\sum_{i=a}^k m_i
\qquad (1\le a\le k), 
\en
and $\lambda(\mb)=(\lambda_1(\mb),\cdots,\lambda_k(\mb))$. 
Let
$\Mb=(M_1,\cdots,M_k), \mb=(m_1,\cdots,m_k)\in\Z_{\ge 0}^k$,  
$\lambda=\lambda(\Mb)$, $\mu=\lambda(\mb)$.
We define
\bea
F_{\Mb,\mb}(q)
&=&
q^{\sum_{a=1}^{k-1}\mu_{a+1}(\lambda_a-\mu_a)}
\prod_{a=1}^{k}
\left[{\lambda_a-\mu_{a+1}\atop\mu_a-\mu_{a+1}}\right]. 
\label{F}
\ena
These are $q$-analogs of the coefficients $F_{\Mb,\mb}(1)$
appearing in the expansion
\bea
&&(1+x_1)^{M_1}(1+x_1+x_2)^{M_2}\cdots(1+x_1+\cdots+x_k)^{M_k}
\nn\\
&&=
\sum_{\mb}
F_{\Mb,\mb}(1)x_1^{m_1}x_2^{m_2}\cdots x_k^{m_k}.
\label{super}
\ena
Here the sum is taken over $\mb
=(m_1,\cdots,m_k)\in\Z_{\ge 0}^k$. 
Note that the summand is zero unless
$\lambda_a(\mb)\le\lambda_a(\Mb)$ for $a=1,\cdots,k$. 

Another ingredient is the level-restricted Kostka polynomial 
for $\slt$.  
We make use of its fermionic formula \cite{SS} 
given by 
\bea
K^{(k)}_{l,\mb}(q)
&=&
\sum_{\nb\in \Z_{\ge 0}^k
\atop 2|\nb|=|\mb|-l}
q^{c(\nb)}\prod_{a=1}^k
\qbin{p_a+n_a}{n_a}.
\label{resKos}
\ena
Here we have set $\nb=(n_1,\cdots,n_k)$ and 
\be
&&c(\nb)=
\sum_{a,b=1}^kA_{ab}n_a n_b
+\sum_{a=1}^kv_an_a,
\\
&&p_a=\sum_{a,b=1}^k A_{ab}(m_b-2n_b)-v_a,
\\
&&
A_{ab}=\min(a,b),
\quad
v_a=\max(a-k+l,0).
\label{GordonA}
\en

The third ingredient is the character of the fusion 
product of irreducible $\slt$-modules. 
Let 
$\pi_l$ denote
the $(l+1)$-dimensional irreducible module of
$\slt$, and let $u_l$ be the lowest weight vector. 
The fusion product 
\bea
\pi_\mb=\pi_1^{*m_1}*\cdots*\pi_k^{*m_k}, 
\label{fusslt}
\ena
where $\pi_a^{*m_a}$ means
$\pi_a*\cdots*\pi_a$ ($m_a$-times), 
is a graded $\slt[t]$-module. 
Assign degrees \eqref{degree} 
and the degree $(0,-|\mb|/2)$ to the cyclic vector 
$u_1^{\otimes m_1}\otimes\cdots\otimes u_k^{\otimes m_k}$. 
The following formula is due to \cite{FL},\cite{FF}.
\bea
&&\ch_{q,z}\pi_\mb
=z^{-|\mb|/2}\chi_{\mb}(q,z),
\nn\\
&&\chi_{\mb}(q,z)
=\sum_{\nb}z^{|\mb|-|\nb|}F_{\mb,\nb}(q). 
\label{chi1}
\ena

\subsection{Dimensions}\label{subsec:2.4}
As mentioned in Introduction, 
the dimensions of the spaces of coinvariants
\eqref{bigc},\eqref{mixc} are described by the use of 
the Verlinde algebra. 
In the case \eqref{bigc} we have 
\bea
\dim \Bigl(L_l^{(k)}/B_N\Bigr)&=&
\bigl((\sum_{j=0}^k(j+1)[j])^N:[l]\bigr)_k
\nn\\
&=&\sum_{\mb}
F_{(0,\cdots,0,N),\mb}(1)
K^{(k)}_{l,\mb}(1)\chi_{\mb}(1,1). 
\label{dimbig}
\ena
Similarly, in the case \eqref{mixc} we have 
\bea
\dim \Bigl(L_l^{(k)}/Y_{\Mb,\bMb}(\zz)\Bigr)
&=&
\bigl(
\prod_{a=1}^k([0]+[1]+\cdots+[a])^{M_a+\bar{M}_a}:[l]
\bigr)_k
\nn\\
&=&\sum_{\mb, \bmb}F_{\Mb,\mb}(1)F_{\bMb,\bmb}(1)
K^{(k)}_{l,\mb+\bmb}(1).
\label{dimmix}
\ena
In both of these formulas, the 
second equality follows from the first by 
using \eqref{super} and the known formula 
\be
\bigl([1]^{m_1}\cdots [k]^{m_k}:[l]\bigr)_k
=K^{(k)}_{l,\mb}(1).
\en

Eq. \eqref{dimbig} was obtained in \cite{FL}. 
Eq. \eqref{dimmix} is a special case of 
a multiplicative formula obtained in \cite{FJKLM1}  
(Theorem 2.9, see also Example 4). 

\subsection{Fermionic formulas}\label{subsec:2.5}
We now state our main results which are
natural $q$-analogs of the above dimension formulas.

\begin{thm}\label{thm:2.1}
The character of the space \eqref{bigc} is given by 
\bea
\ch_{q,z} \bigl(L_l^{(k)}/B_N \bigr)
&=&\sum_{\mb}
F_{(0,\cdots,0,N),\mb}(q)K^{(k)}_{l,\mb}(q)
\ch_{q,z}\pi_\mb.
\label{chbig}
\ena
\end{thm}

\begin{thm}\label{thm:2.2}
The character of the associated graded space of 
\eqref{mixc} is given by
\bea
&&\ch_{q,z}\ggr\left(L_l^{(k)}/Y_{\Mb,\bMb}(\zz)
\right)
=\sum_{\mb,\bmb}F_{\Mb,\mb}(q)F_{\bMb,\bmb}(q)
K^{(k)}_{l,\mb+\bmb}(q)
z^{|\mb|-|\bmb|}.
\nn\\
&&
\label{chmix}
\ena
\end{thm}

In the previous studies, 
fermionic formulas of different type have been obtained:   
for \eqref{bigc} in \cite{FKLMM1} (when $k=1$),   
and for \eqref{mixc} in \cite{FKLMM3}
(when $M_{ik}=\delta_{ik}M,\bM_{ik}=\delta_{ik}\bM$). 
Theorem \ref{thm:2.2} is a direct generalization of 
the formula (3.15) in \cite{FJKLM2} 
where the case \eqref{mixc} with 
$\bM_1=\cdots=\bM_k=0$ was treated.  

The rest of the text is devoted to the 
proof of Theorems \ref{thm:2.1}, \ref{thm:2.2}. 

\section{Space of coinvariants $L_l^{(k)}/B_N$}\label{sec:3}
In this section we prove Theorem \ref{thm:2.1}. 
\subsection{Coinvariants of fusion products}\label{subsec:3.1}
{}For the computation of the character \eqref{chbig}, 
we make use of 
the fusion product of certain reducible $\slt$-modules.  

{}Fix $N$ distinct complex numbers 
$\zz=(\z_1,\cdots,\z_N)$.
Set 
\be
B_{1,\zz}=\ga_N(\zz)U(\slt[t]),
\quad    
\ga_N(\zz)=
\slt\otimes\prod_{j=1}^N(t-\z_j)\C[t].
\en
Since $B_N$ and $B_{1,\zz}$ are two-sided ideals, 
the space of coinvariants $L_l^{(k)}/B_N$,
$L_l^{(k)}/B_{1,\zz}$ admit the action of $\slt[t]$. 
The former may be viewed as the `limit' of 
the latter when all points $\z_i$ tend to $0$. 
It is known that the dimension does not change in this limit: 

\begin{thm}\label{thm:3.1}$($\cite{FKLMM1}, Theorem 9$)$
We have the equality of dimensions 
\be
\dim\bigl(L_l^{(k)}/B_N\bigr)=
\dim\bigl(L_l^{(k)}/B_{1,\zz}\bigr).
\en
\end{thm}

Let $\pi_l^\vee$ denote the dual representation of 
$\pi_l$, on which $\slt$ acts from the left.  
Consider a reducible $\slt\oplus\slt$-module 
\bea
\vpi^{(k)}=\bigoplus_{l=0}^k\pi_l^\vee\otimes\pi_l. 
\label{pik}
\ena
Let $\su_l\in \pi_l^\vee\otimes\pi_l$ be the canonical vector.
Then   
\bea
\su^{(k)}=\sum_{l=0}^k \su_l~\in \vpi^{(k)}
\label{scyc}
\ena
is a cyclic vector of \eqref{pik}
viewed as a module over $0\oplus\slt$. 
We consider the filtered tensor product 
$\F_\zz(\vpi^{(k)},\cdots,\vpi^{(k)})$ and the associated 
fusion product  
\bea
\vpi^{(k)}*\cdots*\vpi^{(k)}, 
\label{fus0}
\ena
by choosing 
$\sub=\su^{(k)}\otimes\cdots\otimes \su^{(k)}$
as the cyclic vector.  
{}For $x\in\slt$, write
$x'=(x,0), x''=(0,x)\in\slt\oplus\slt$. 
Since $(x'+x'')\sub=0$, the filtered tensor
product as $(\slt\oplus\slt)[t]$-modules 
coincides with the one as $(0\oplus\slt)[t]$-modules.  
\newpage

\begin{thm}\label{thm:3.2}
$($\cite{FJKLM1},Theorem 3.6, Theorem A.3$)$
\begin{enumerate}
\item
We have an isomorphism of filtered $\slt[t]$-modules
\be
L_l^{(k)}/B_{1,\zz}\simeq
\F_\zz(\vpi^{(k)},\cdots,\vpi^{(k)})/
\langle e''_0,{e_1''}^{k-l+1},h''_0+l\rangle,
\en
where $\langle S\rangle$
signifies the right ideal of $U(0\oplus\slt[t])$
generated by the set $S$. 
The action of an element $x_i\in \slt[t]$ on the left hand side 
is sent to the action of $x_i'$ on the right hand side.
\item There are canonical surjections
\be
L_l^{(k)}/B_N 
\longrightarrow
\vpi^{(k)}*\cdots*\vpi^{(k)}/
\langle e''_0,{e''_1}^{k-l+1},h''_0+l\rangle 
\longrightarrow
\gr \Bigl(L_l^{(k)}/B_{1,\zz}\Bigr).  
\en
\end{enumerate}
\end{thm}

Combining Theorem \ref{thm:3.1} and
Theorem \ref{thm:3.2} we obtain 

\begin{prop}\label{prop:3.1}
We have an isomorphism of graded $\slt[t]$-modules
\be
L_l^{(k)}/B_N\simeq \vpi^{(k)}*\cdots*\vpi^{(k)}/
\langle e''_0,{e''_1}^{k-l+1},h''_0+l\rangle.
\en
In particular, both sides have the same characters. 
\end{prop}
In what follows we study the fusion product \eqref{fus0} 
and its quotient. 

\subsection{Changing cyclic vectors}\label{subsec:3.2}
As the next step, 
we change the cyclic vector of \eqref{fus0} 
into a simpler one. 
{}For that purpose we utilize the embedding of Lie algebras 
$\slt\oplus\slt\simeq\go_4\subset\go_5$. 

The Lie algebra $\go_5$ 
is realized as the Lie algebra of matrices
\be
\go_5=
\{\sum_{a,b=1}^5c_{ab}e_{ab}\mid c_{ab}\in\C,~~
c_{ab}+c_{6-b\,6-a}=0\quad(1\le a,b\le 5)\}.
\en
The set of positive roots of $\go_5$ has the form 
$\Delta_+=\{\ep_1-\ep_2,\ep_1+\ep_2,\ep_1,\ep_2\}$, 
where $\ep_1,\ep_2$ are orthonormal vectors. 
We choose root vectors 
corresponding to $\Delta_+$ and $-\Delta_+$ as follows.
\be
&X=e_{12}-e_{45},\quad
&\overline{X}=e_{21}-e_{54},
\\
&Y=e_{14}-e_{25},\quad
&\overline{Y}=e_{41}-e_{52},
\\
&Z=\sqrt{2}(-e_{13}+e_{35}),\quad
&\overline{Z}=\sqrt{2}(-e_{31}+e_{53}),
\\
&T=\sqrt{2}(e_{23}-e_{34}), \quad
&\overline{T}=\sqrt{2}(e_{32}-e_{43}). 
\en
The abelian subalgebra
\be
\ga=\C X\oplus\C Y\oplus \C Z\,\subset\, \go_5
\en
will play a role in the sequel.
We have a commuting pair of $\slt$
\be
&&\slt'=\C Y\oplus\C \bY\oplus \C [Y,\bY],
\\
&&\slt''=\C X\oplus\C \bX\oplus \C [X,\bX],
\en
which span the subalgebra $\go_4\subset \go_5$. 

Let $\omega_1=\ep_1$ be the highest weight of 
the natural representation $\C^5$ of $\go_5$.  
Denote by $\pit^{(k)}$ the representation with 
highest weight $k\omega_1$, and let $\ut^{(k)}$ 
be its lowest weight vector. 
Denoting by
$\overline{\gn}$ the span of the negative root vectors 
$\overline{X},\overline{Y},\overline{Z},\overline{T}$, 
we have 
\bea
&&\overline{\gn}\,\ut^{(k)}=0,
\label{ann1}\\
&&
T \ut^{(k)}=0, \quad X^{k+1} \ut^{(k)}=0.
\label{ann2}
\ena

\begin{lem}\label{lem:3.2}
\begin{enumerate}
\item 
As a module over the subalgebra $\go_4$ we have a decomposition
\bea
\pit^{(k)}\simeq\bigoplus_{l=0}^k\pi_l\otimes\pi_l.
\label{decomp1}
\ena
\item 
Regard $\Pi^{(k)}$ as $\slt$-module via the action of 
$\slt''$, and $\varpi^{(k)}$ via that of $0\oplus\slt$.  
Then there exists an isomorphism of $\slt$-modules 
\be
\nu:\Pi^{(k)} \overset{\sim}{\longrightarrow}
\varpi^{(k)}=\bigoplus_{l=0}^k\pi_l^{\vee}\otimes\pi_l
\en
such that $\nu(e^Z \ut^{(k)})=\sigma^{(k)}$. 
\end{enumerate}
\end{lem}
\begin{proof}
Assertion (i) can be easily verified, say by comparing 
characters. 

To see (ii), first note the following. 
Let $\theta$ be the anti-automorphism of $\slt$ given by
$\theta(e)=f,\theta(h)=h$. 
{}Fix a non-degenerate symmetric bilinear form $(~|~)$ on 
$\pi_l$ satisfying $(\theta(x)u|v)=(u|xv)$ for 
$x\in\slt$, $u,v\in \pi_l$, 
and define $i:\pi_l\rightarrow\pi_l^{\vee}$ by 
$i(u)=(u|\cdot)$. 
Then an element $\sigma\in \pi_l^{\vee}\otimes\pi_l$
is proportional to the canonical element if and only if
the element
$\sigma'=(i^{-1}\otimes 1)\sigma\in \pi_l\otimes\pi_l$ 
satisfies $(-\theta(x)\otimes 1+1\otimes x)\sigma'=0$
for $x\in\slt$. 

We have
\be
&&(-\overline{Y}+X)e^Z \ut^{(k)}=
e^Z (-\overline{Y}-\overline{T})\ut^{(k)}=0,
\\
&&(-Y+\overline{X})e^Z \ut^{(k)}=
e^Z (-T+\overline{X})\ut^{(k)}=0.
\\
\en
We will show also 
in the proof of Proposition 
\ref{prop:3.2} below that 
$e^Z \ut^{(k)}$ generates $\Pi^{(k)}$ over $\slt''$. 
Hence if $e^Z \ut^{(k)}=\sum_{l=0}^k\sigma'_l$ 
is the decomposition according to \eqref{decomp1}, then 
we have $\sigma_l'\neq 0$. 
Assertion (ii) follows from these observations. 
\end{proof}

Let now $k_i\in \{0,1,\cdots,k\}$, $i=1,\cdots,N$. 
Generalizing slightly the setting of \eqref{fus0}, 
we consider the filtered tensor product of $\go_5[t]$-modules
\bea
\F_\zz\left(\pit^{(k_1)},\cdots,\pit^{(k_N)}\right),
\label{filt9}
\ena
choosing the cyclic vector 
\bea
\vb=\ut^{(k_1)}\otimes\cdots\otimes\ut^{(k_N)}.
\label{cycvb}
\ena
Set
\be
\sub=\su^{(k_1)}\otimes\cdots\otimes \su^{(k_N)}.  
\en
We retain the notation $Z_i=Z\otimes t^i$ and so forth. 
\newpage

\begin{prop}\label{prop:3.2}
The following are isomorphic as filtered vector spaces.
\begin{enumerate}
\item Filtered tensor product of $(0\oplus\slt)[t]$-modules 
$\F_\zz(\vpi^{(k_1)},\cdots,\vpi^{(k_N)})$, with 
$\sub$ as cyclic vector,   
\item Filtered tensor product of $\go_5[t]$-modules 
$\F_\zz(\pit^{(k_1)},\cdots,\pit^{(k_N)})$, with 
$\vb$ as cyclic vector,   
\item Filtered tensor product of $\ga[t]$-modules 
$\F_\zz(\pit^{(k_1)},\cdots,\pit^{(k_N)})$, 
with $\vb$ as cyclic vector.    
\end{enumerate}
\end{prop}

{}For the proof we use 
\begin{lem}\label{lem:3.3}
{}For all $m$ and $i\ge 0$ we have
\be
&&T_iZ_0^{(m)}\vb 
=-2Y_iZ_0^{(m-1)}\vb,
\\
&&\bX_i Z_0^{(m)}\vb
=Y_iZ_0^{(m-2)}\vb, 
\en
where $Z_0^{(m)}=Z_0^m/m!$ if $m\ge 0$ and $=0$ if $m<0$. 
\end{lem}
\begin{proof}
This follows from $[\bX_i,Z_0]=-T_i$, 
$[T_i,Z_0]=-2Y_i$, $\bX_i\vb=T_i\vb=0$
and that $\bX_i,T_j,Y_l$ are mutually commutative.
\end{proof}

\noindent{\it Proof of Proposition \ref{prop:3.2}.}\quad
The equivalence of (ii) and (iii) is a consequence of 
the relations \eqref{ann1},\eqref{ann2}
and the Poincar{\' e}-Birkhoff-Witt (PBW) theorem. 

{}For a Lie algebra $\g$, let 
$U^{\le d}(\g[t])$ stand for the subspace of 
$U(\g[t])$ spanned by elements of degree at most $d$  
with respect to the grading in $t$. 
Set $F^{d}=U^{\le d}(\slt''[t])e^{Z_0}\vb$,  
$G^{d}=U^{\le d}(\go_5[t])\vb$.  
In view of Lemma \ref{lem:3.2}, 
to show the equivalence of (i) and (ii) 
it suffices to prove 
the equality $F^{d}=G^{d}$ for all $d\ge 0$.  
Since $e^{Z_0}\vb\in G^0$ 
and $\slt''[t]\subset\go_5[t]$, 
we have $F^d\subset G^d$. Let us prove the opposite inclusion.

Set $H_p=[X_p,\bX_0]\in \slt''[t]$. 
Since $[H_0,Z_0^{(m)}]=mZ_0^{(m)}$,  
$e^{Z_0}\vb\in F^0$ implies
$Z_0^{(m)}\vb\in F^0$ for all $m$. 
{}From Lemma \ref{lem:3.3} we have
\be
&&Y_iZ_0^{(m)}\vb=\bX_iZ_0^{(m+2)}\vb,
\\
&&T_iZ_0^{(m)}\vb=-2\bX_iZ_0^{(m+1)}\vb.
\en
Since $X_i F^d, \bX_i F^d\subset F^{d+i}$ and $\bX_i,T_j,Y_l$ 
are mutually commutative, we find by induction that 
\be
\Bigl(\prod_{a=1}^rX_{i_a}\prod_{b=1}^s Y_{j_a}\prod_{c=1}^t T_{l_c}\Bigr)
Z_0^{(m)}\vb\in F^d
\en
if $\sum_{a=1}^ri_a+\sum_{b=1}^sj_b+\sum_{c=1}^tl_c\le d$. 
Applying $H_p$ repeatedly to this expression
using $[H_p,Z_i]=Z_{i+p}$, we 
obtain 
\be
\Bigl(\prod_{a=1}^rX_{i_a}\prod_{b=1}^sY_{j_a}\prod_{c=1}^t
T_{l_c}
\prod_{g=1}^wZ_{p_g}\Bigr)\vb\in F^d
\en
for 
$\sum_{a=1}^ri_a+\sum_{b=1}^sj_b+\sum_{c=1}^tl_c+
\sum_{g=1}^wp_g\le d$. 
Therefore the PBW theorem implies $G^d\subset F^d$.
\qed

Henceforth we set 
\bea
&&V_\Mb=\pit^{(k_1)}*\cdots*\pit^{(k_N)},
\label{fusV}
\\
&&
\Mb=(M_1,\cdots,M_k), \quad
M_a=\sharp\{j\mid k_j=a\}. 
\label{Mb}
\ena
We regard the fusion product as a module over  
$U(\ga[t])=\C[X_i,Y_i,Z_i~~(i\ge 0)]$
with cyclic vector \eqref{cycvb}.
{}From \eqref{super} and \eqref{decomp1} we have
\bea
\dim V_\Mb
&=&\dim \Pi^{(k_1)}\otimes\cdots\otimes \Pi^{(k_N)}
\nn\\
&=&
\sum_{\mb}
F_{\Mb,\mb}(1)\left(\dim \pi_1^{\otimes m_1}\otimes
\cdots\otimes\pi_k^{\otimes m_k}\right)^2.
\label{dimV}
\ena
The operators $D_0=td/dt$, $D_1=[Y_0,\overline{Y}_0]/2$, 
$D_2=[X_0,\overline{X}_0]/2$ give the following
grading on \eqref{fusV}.
\bea
&&\deg X_i=(i,0,1),
\\
&&\deg Y_i=(i,1,0),
\\
&&\deg Z_i=(i,\frac{1}{2},\frac{1}{2}),
\\
&&\deg \vb=(0,-\frac{|\Mb|}{2},-\frac{|\Mb|}{2}).
\ena 
Our goal is to obtain the character of the quotient space
\bea
V_\Mb/(X_0V_\Mb+X_1^{k-l+1}V_\Mb+(D_2+l/2)V_\Mb).
\label{goal}
\ena

\subsection{Annihilation conditions}\label{subsec:3.3}
Let us determine the relations satisfied by the
cyclic vector $\vb\in V_\Mb$. 

{}From the relations \eqref{ann1}, \eqref{ann2}, 
for any $\xi$ we have 
$e^{\xi T}X^{k+1}e^{-\xi T}v^{(k)}=0$, or equivalently
\be
\left(X+\xi Z-\xi^{2}Y\right)^{k+1}v^{(k)}=0.
\en

{}For an element $\eta\in\ga$, 
consider the generating function
$\eta(z)=\sum_{i\ge 0}\eta_i z^i$.
On the $N$-fold fusion product \eqref{fusV}, 
the operators $\eta_i$ ($i\ge N$) act as $0$. 
Therefore $\eta(z)$ acts as
a polynomial of degree at most $N-1$. 

\begin{prop}\label{prop:3.3}
{}For an indeterminate $\xi$, 
the following relations hold on the fusion product $V_\Mb$:
\bea
&&
\deg_z\bigl(X(z)+\xi Z(z)-\xi^2Y(z)\bigr)^\nu\vb
\le\sum_{i=1}^k\min(\nu,i)M_i-\nu
\label{acond1}
\\
&&\qquad \mbox{for any $\nu\ge 0$. }\nn
\ena
Here $\deg_z$ signifies the degree of a polynomial in $z$.
\end{prop}
\begin{proof}
We repeat the argument of \cite{FF}. 
{}For convenience we assume that $\z_a\neq0$. 
Set $x(z)=X(z)+\xi Z(z)-\xi^2Y(z)$ and 
$\xt(z)=\prod_{a=1}^N(1-\z_az)\cdot x(z)=\sum_{i\ge 0}\xt_iz^i$.
By the definition of the filtered tensor product,  
$x_i$ acts on the $a$-th tensor component of \eqref{filt9}
as $\z_a^i x$. 
We have an operator identity 
\be
\xt(z)=\sum_{a=1}^N\prod_{b(\neq a)}(1-\z_{b}z)\,
1\otimes\cdots \otimes 
\overset{\mbox{\tiny $a$-th}}{x}
\otimes\cdots\otimes 1 .
\en
In particular, the relations 
$\xt_i=0$ ($i\ge N$) and $\xt(\z_a^{-1})^{k_a+1}=0$ hold.  
Therefore $\xt(z)^\nu$ is divisible by
$(1-\z_az)^{\max(\nu-k_a,0)}$, 
so that
\be
\prod_{a=1}^N(1-\z_az)^{\nu-\max(\nu-k_a,0)}\cdot
x(z)^\nu
\en
acts as a polynomial of degree at most
\be
\nu (N-1)-\sum_{a=1}^N \max(\nu-k_a,0)
=\sum_{a=1}^N\min(\nu,k_a)-\nu.
\en
Passing to the associated graded space, 
$1-\z_az$ can be replaced by $1$ because
$z$ picks up operators of lower degree. 
The assertion follows from this.
\end{proof}

The relations \eqref{acond1} are equivalently written in the form
\be
&&\deg_z\Bigl(X(z)^aZ(z)^b+\sum_{1\le j\le b/2}
C^{ab}_jX(z)^{a+j}Y(z)^jZ(z)^{b-2j}\Bigr)\vb
\\
&&
\le\sum_{i=1}^k\min(a+b,i)M_i-(a+b)
\qquad (a,b\ge 0)
\en
with some constants $C^{ab}_j$, 
and similar relations with $X$ and $Y$ interchanged.
Let us simplify them further.
Denote by $\cG^m$ the $\C[X_i,Y_i~~(i\ge 0)]$-submodule
of $V_\Mb$ spanned by 
$\prod_{a=1}^rX_{i_a}\prod_{b=1}^s Y_{j_b}
\prod_{c=1}^t Z_{l_c}\vb$ with $t\le m$. 
This defines a filtration 
\be
\cG~:~
0=\cG^{-1}\subset \cG^0\subset \cdots \subset \cG^{i-1}
\subset \cG^{i}\subset \cdots\subset V_\Mb.
\en
On the associated graded space 
$\ggr^\cG V_\Mb=\oplus_{i\ge 0} \cG^i/\cG^{i-1}$,  
we have the relations 
\bea
&&\deg_z X(z)^aZ(z)^b\vb
\le\sum_{i=1}^k\min(a+b,i)M_i-(a+b),
\label{acond2}\\
&&\deg_z Y(z)^aZ(z)^b\vb
\le\sum_{i=1}^k\min(a+b,i)M_i-(a+b),
\label{acond3}
\ena
valid for all $a,b\ge 0$.

\subsection{Subquotient modules and recursion}\label{subsec:3.4}
Suggested by the relations \eqref{acond2},\eqref{acond3}, 
we introduce a family of cyclic modules 
\bea
W(k_1,\cdots,k_p|l_1,\cdots,l_r)
=
\C[X_i,Y_i,Z_i~~(i\ge 0)]\bone, 
\label{W0}
\ena
defined by the following relations for the cyclic vector
$\bone$:  
\bea
&&\deg_z X(z)^aZ(z)^b \bone
\le \sum_{i=1}^p\min(a+b,k_i)+
\sum_{j=1}^r\min(a,l_j)-(a+b),
\label{acond4}\\
&&\deg_z Y(z)^aZ(z)^b \bone
\le \sum_{i=1}^p\min(a+b,k_i)+
\sum_{j=1}^r\min(a,l_j)
-(a+b),
\label{acond5}
\ena
for all $a,b\ge 0$.
We also write \eqref{W0} as 
\be
&&W_{\Mb,\nb}=W(k_1,\cdots,k_p|l_1,\cdots,l_r),
\\
&&\nb=(n_1,\cdots,n_k),\quad n_a=\sharp\{j\mid l_j=a\}. 
\en

In the case $p=N$ and $r=0$, 
\eqref{acond2},\eqref{acond3} imply that   
we have a surjection 
\bea
W(k_1,\cdots,k_N|)\longrightarrow
\ggr^\cG V_M\longrightarrow 0.
\label{surj}
\ena
In the case $p=0$, the relations 
\eqref{acond4},\eqref{acond5} reduce to 
\be
&&Z(z)\bone=0,
\\
&&\deg_z X(z)^a\bone\le 
\sum_{j=1}^r\min(a,l_j)-a,
\\
&&\deg_z Y(z)^a\bone\le 
\sum_{j=1}^r\min(a,l_j)-a.
\en
The last two relations are each identical to the defining relations 
for the cyclic vector of 
the fusion product of irreducible $\slt[t]$-modules \cite{FF}.  
Therefore we have
\bea
W(|l_1,\cdots,l_r)
\simeq
(\pi_{l_1}*\cdots*\pi_{l_r})
\otimes
(\pi_{l_1}*\cdots*\pi_{l_r})
\label{ini}
\ena
as $\C[X_i,Y_i~~(i\ge 0)]$-modules.  

We are now in a position to give a recursion relation among
\eqref{W0}. 
Suppose $p\ge 1,k_1\ge 1$, and 
set
\be
&&
W=W(k_1,\cdots,k_p|l_1,\cdots,l_r),
\\
&&
W'=W(k_2,\cdots,k_p|k_1,l_1,\cdots,l_r),
\\
&&
W''=W(k_1-1,k_2,\cdots,k_p|l_1,\cdots,l_r).
\en
Let further $\bW'$ be the subspace 
of $W$ generated by $\bone$ over 
$X_i,Y_i$ ($i\ge 0$) and $Z_i$ ($i\ge 1$).

The following Proposition can be shown in exactly the same 
way as that of Propositions 2.6 and 2.8 in \cite{FJKLM2}.  
\begin{prop}\label{prop:3.4}
\begin{enumerate}
\item There exists a surjection 
\be
\iota~:~W'\longrightarrow \bW'
\en
given by $\iota(X_i)=X_i$, $\iota(Y_i)=Y_i$, $\iota(Z_i)=Z_{i+1}$.
\item Assume that $k_1\ge k_2,\cdots,k_N$. 
Then there exists a surjection 
\be
\phi~:~W''\longrightarrow W/\bW'
\en
given by $\phi(bZ_0^{(m)})=b Z_0^{(m+1)}$, 
where $b$ is an element not divisible by $Z_0$. 
\end{enumerate}
\end{prop}

Introduce a $\Z^4_{\ge 0}$ grading on \eqref{W0} 
\bea
&&\deg X_i=(i,0,1,0),
\quad
\deg Y_i=(i,1,0,0),
\quad
\deg Z_i=(i,0,0,1),
\label{z4}\\
&&\deg {\bf 1}=(0,0,0,0).
\label{z42}
\ena
Denoting by $W_{d,i_1,i_2,i_3}$ the homogeneous component of
degree $(d,i_1,i_2,i_3)$ we have 
\be
&&\iota(W'_{d-i_3,i_1,i_2,i_3})\subset W_{d,i_1,i_2,i_3},
\\  
&&\phi(W''_{d,i_1,i_2,i_3-1})\subset (W/\bW')_{d,i_1,i_2,i_3}.
\en

\begin{prop}\label{prop:3.41}
We have an upper estimate
\bea
&&\ch_{q,z_1,z_2,z_3} W_{\Mb,\nb}
\nn\\
&&\le
\sum_{\mb}
 F_{\Mb,\mb}(q)z_3^{|\Mb|-|\mb|}
\chi_{\mb+\nb}(q,z_1)\chi_{\mb+\nb}(q,z_2), 
\label{esti}
\ena
where $\chi_{\mb}(q,z)$ is given in \eqref{chi1}.
\end{prop}

\begin{proof}
Suppose $k_1\ge k_2,\cdots,k_p$. 
{}From Proposition \ref{prop:3.4} we have an exact sequence 
of $\C[X_i,Y_i~~(i\ge 0)]$-modules
\be
&&
\begin{CD}
\phantom{0}@.W'@.@.W''@.
\\
@.@V{\iota}VV @.\phantom{W}@VV{\phi}V @.@.
\\
0@>>>\overline{W}' @>>> W @>>> W/\overline{W}' @>>>0
\end{CD}
\en
where the vertical arrows are surjective. 
Hence we have 
\be
&&\ch_{q,z_1,z_2,z_3} W(k_1,\cdots,k_p|l_1,\cdots,l_r)
\\
&&\le
\ch_{q,z_1,z_2,qz_3}W(k_2,\cdots,k_p|k_1,l_1,\cdots,l_r)
\\
&&+
z_3\,\ch_{q,z_1,z_2,z_3} W(k_1-1,k_2,\cdots,k_p|l_1,\cdots,l_r). 
\en
Repeating the working of Theorem 2.11 
in \cite{FJKLM2} and using 
\be
\ch_{q,z_1,z_2,z_3}W(|l_1,\cdots,l_r)
=\chi_{\nb}(q,z_1)\chi_{\nb}(q,z_2)
\en
which follows from \eqref{ini}, 
we obtain the assertion.
\end{proof}

In Proposition \ref{prop:3.41}, take $r=0$ and 
specialize to $q=z_1=z_2=z_3=1$. 
Using \eqref{dimV}, we obtain an estimate 
\be
\dim W(k_1,\cdots,k_N|)\le \dim V_\Mb. 
\en
{}From \eqref{surj} we have also the opposite inequality. 
We thus find that 
\be
W(k_1,\cdots,k_N|)\simeq \ggr^\cG V_\Mb.
\en 
At the same time, the maps $\iota,\phi$ appearing in the
intermediate steps are isomorphisms. 
This implies that $W(k_1,\cdots,k_N|)$ has a filtration 
with subquotients of the form $W(|l_1,\cdots,l_r)$. 

Choosing $z_3^2=z_1z_2$ in 
the right hand side of \eqref{esti}, 
and multiplying $(z_1z_2)^{-|\Mb|/2}$ on both sides, 
we obtain a formula for the character of $V_\Mb$.  
 
Let us summarize the conclusion as Theorem. 
\begin{thm}\label{thm:3.3}
There exists a filtration $\cH$ of 
$\ggr^\cG V_\Mb$ by $\C[X_i,Y_i~~(i\ge 0)]$-modules 
such that 
\be
&&\ggr^{\cH}\ggr^\cG V_\Mb
=\bigoplus_{\mb}\mathcal{M}_{\Mb,\mb}
\otimes \pi_\mb\otimes\pi_\mb,
\en
where $\mathcal{M}_{\Mb,\mb}$ is a trivial module with
the character
\be
&&\ch_q\mathcal{M}_{\Mb,\mb}=F_{\Mb,\mb}(q).
\en
In particular, we have 
\be
\ch_{q,z_1,z_2}V_{\Mb}
=\sum_{\mb}
F_{\Mb,\mb}(q)\,
\ch_{q,z_1}\pi_{\mb}\,\ch_{q,z_2}\pi_{\mb}.
\en
\end{thm}

\subsection{Proof of Theorem \ref{thm:2.1}}\label{subsec:3.5}
Let us return to the space of coinvariants \eqref{bigc}.  
\medskip

\noindent{\it Proof of Theorem \ref{thm:2.1}.}\quad
It remains to calculate the character of 
the quotient space \eqref{goal}.
Let us set $\Mb=(0,\cdots,0,N)$, 
$V=V_\Mb$, $V'=\ggr^\cG V$ and  
$V''=\ggr^\cH\ggr^\cG V$. 
We have
\be
&&\ch_{q,z}V/(X_0V+X_1^{k-l+1}V+(D_2+l/2)V)
\\
&&\le 
\ch_{q,z}V'/(X_0V'+X_1^{k-l+1}V'+(D_2+l/2)V')
\\
&&\le 
\ch_{q,z}V''/(X_0V''+X_1^{k-l+1}V''+(D_2+l/2)V'')
\\
&&
=\sum_{\mb}
F_{\Mb,\mb}(q)\,
\ch_{q,z}\pi_{\mb}
\\
&&\qquad \times 
\ch_{q}
\pi_{\mb}/(X_0 \pi_\mb+X_1^{k-l+1}\pi_\mb+(D_2+l/2)\pi_\mb).
\en
The last equality follows from Theorem \ref{thm:3.3}.
Using the formula (\cite{FJKLM1}, Theorem 4.1) 
\be
\ch_{q}\pi_{\mb}
/(e_0\pi_{\mb}+e_1^{k-l+1}\pi_{\mb}+(h_0+l)\pi_{\mb})
=K_{l,\mb}^{(k)}(q), 
\en
we find that the last line coincides with 
the right hand side of \eqref{chbig}.  
Comparing dimensions, we obtain the desired equality.
Theorem \ref{thm:2.1} is proved. 
\qed

Let us mention an immediate consequence of the above proof.  
Set 
\be
\cK_{l,\Mb}(q,z)=\ch_{q,z}V_{\Mb}
/(e_0V_{\Mb}+(h_0+l)V_{\Mb}).
\en

\begin{cor}\label{cor:conj}
{}For $\Mb=(0,\cdots,0,N)$, we have
\be
\ch_{q,z} \bigl(L_l^{(k)}/B_N \bigr)
&=&
\sum_{i\ge 0}q^{(k+2)i^2+(l+1)i}
\cK_{2(k+2)i+l,\Mb}(q,z) 
\\
&&-\sum_{i> 0}q^{(k+2)i^2-(l+1)i}
\cK_{2(k+2)i-l-2,\Mb}(q,z).  
\en
\end{cor}
\begin{proof}
Let 
$K_{l,\mb}(q)=\ch_q\pi_\mb/(e_0\pi_\mb+(h_0+l)\pi_\mb)$
denote the (non-restricted) Kostka polynomial.
The following alternating sum formula is known
(\cite{SS}, eq.(6.8)):
\be
K_{l,\mathbf{m}}^{(k)}(q) 
&=&
\sum_{i\ge 0}q^{(k+2)i^2+(l+1)i}
K_{2(k+2)i+l,\mathbf{m}}(q) 
\\
&&-\sum_{i> 0}q^{(k+2)i^2-(l+1)i}
K_{2(k+2)i-l-2,\mathbf{m}}(q).
\en
Substituting this into \eqref{chbig} we obtain the assertion.
\end{proof}

\begin{rem}
Corollary \ref{cor:conj} confirms a conjecture of 
\cite{FJKLM1}, eq.(3.26), in the special case \eqref{bigc}.  
We remark that a similar alternating sum formula was proposed 
earlier in \cite{FL}. 
While the cyclic vector for \eqref{cor:conj} is the sum of canonical vectors, 
the one in \cite{FL} $($for $\slt$$)$ is chosen to be 
the tensor product of highest weight vectors of $\slt\oplus\slt$.
At this moment we do not know the relation between the two. 
\end{rem}

\section{Space of coinvariants $L_l^{(k)}/Y_{\Mb,\bMb}(\zz)$}\label{sec:4}

The proof of Theorem \ref{thm:2.2} 
is quite parallel to the previous one. 
We describe the main steps below, skipping minor details. 

\subsection{Fusion product}\label{subsec:4.1}
In the case of the space of coinvariants 
\eqref{mixc}, we use reducible modules over
$\C \htil\oplus\slt$, where $\C\htil$ is a one-dimensional 
Lie algebra. 
Let $u_l$, $\bu_l$ denote respectively 
the lowest and highest weight vectors of $\pi_l$.   
Set 
\be
&
\displaystyle{\pi^{(m)}=\bigoplus_{l=0}^m\pi_l}, 
\quad
&
\bpi^{(m)}=\bigoplus_{l=0}^m\pi_l, 
\\
&\displaystyle{u^{(m)}=\sum_{l=0}^mu_l},
\quad
&\bu^{(m)}=\sum_{l=0}^m \bu_l. 
\en
We define the action of $\htil$ on $\pi^{(m)},\bpi^{(m)}$ 
by the rules $\htil u_l=lu_l$, $\htil \bu_l=-l\bu_l$, 
and $[\htil,\slt]=0$.  
Consider the filtered tensor product 
\be
\F_\zz(\pi^{(k_1)},\cdots,\pi^{(k_p)},
\bpi^{(\bk_1)},\cdots,\bpi^{(\bk_{\bp})})
\en 
as $\slt[t]$-modules  
taking $\ub\otimes\bub$ as cyclic vector, where
\bea
\ub=u^{(k_1)}\otimes\cdots \otimes u^{(k_p)},\quad
\bub=\bu^{(\bk_1)}\otimes\cdots \otimes \bu^{(\bk_{\bp})}.
\label{cycv2}
\ena
We use $\Mb,\bMb$ in \eqref{Mb1}, \eqref{Mb2} to label 
the $k_i,\bk_i$. 

\begin{thm}\label{thm:4.1}$($[FJKLM1], Theorem 3.6 and Example 4$)$
We have an isomorphism of filtered vector spaces
\bea
&&L_l^{(k)}/Y_{\Mb,\bMb}(\zz)
\label{Ybb}\\
&&\simeq
\F_\zz(\pi^{(k_1)},\cdots,\pi^{(k_p)},
\bpi^{(\bk_1)},\cdots,\bpi^{(\bk_{\bp})})/
\langle e_0,h_0+l,e_1^{k-l+1}\rangle.
\nn
\ena
The action of $h_0\in\slt[t]$ on the left hand side 
corresponds to that of $\htil$ on the right hand side.
\end{thm}

\subsection{Changing cyclic vectors}\label{subsec:4.2}
We change the cyclic vector \eqref{cycv2} 
to a simpler one with the aid of the embedding $\slt\subset\sltr$.  

Let $\omega_1=\bep_1$, 
$\omega_2=\bep_1+\bep_2=-\bep_3$ 
be the fundamental weights of $\sltr$, 
where $\bep_i=\epsilon_i-(\epsilon_1+\epsilon_2+\epsilon_3)/3$ 
and $\varepsilon_i$ are orthonormal vectors. 
Denote by $\Pi_m$ (resp. $\bPi_m$) the irreducible module with 
highest weight $m\omega_1$ (resp. $m\omega_2$),  
and by $v^{(m)}\in \Pi_m$ (resp. $\bv^{(m)}\in\bPi_m$)
the lowest weight vector. We have 
\be
&&
e_{ij}v^{(m)}=0,\quad 
e_{ij}\bv^{(m)}=0\qquad (i>j),
\\
&&e_{12}v^{(m)}=0,\quad
e_{13}^ae_{23}^bv^{(m)}=0\quad (a+b=m+1),
\\
&&
e_{23}\bv^{(m)}=0,\quad
e_{13}^ae_{12}^b \bv^{(m)}=0\quad (a+b=m+1).
\en
We set 
\bea
&&
\vb=v^{(k_1)}\otimes\cdots\otimes v^{(k_p)},
\quad
\bvb=\bv^{(\bk_1)}\otimes\cdots\otimes \bv^{(\bk_p)}.
\label{cycv3}
\ena
Changing the convention of \cite{FJKLM2}, 
we regard $\slt$ as the subalgebra 
$\C e_{13}\oplus \C e_{31}\oplus\C [e_{13},e_{31}]$ 
of $\sltr$. 
We also use the subalgebras
\be
&&\gn=\C e_{12}\oplus \C e_{13}\oplus \C e_{23}\subset \sltr,
\\
&&\gb=\C e_{13}\subset \gn.
\en

\begin{prop}\label{prop:4.1}
The following are isomorphic as filtered vector spaces.
\begin{enumerate}
\item Filtered tensor product 
of $\slt[t]$-modules  
$\F_{\zz}\left(
\pi^{(k_1)},\cdots,\pi^{(k_p)},
\bpi^{(\bk_1)},\cdots,\bpi^{(\bk_{\bp})}\right)$,
with $\ub\otimes\bub$ as cyclic vector, 
\item Filtered tensor product 
of $\sltr[t]$-modules 
$\F_{\zz}\left(
\Pi_{k_1},\cdots,\Pi_{k_p},
\bPi_{\bk_1},\cdots,\bPi_{\bk_{\bp}}\right)$,  
with $\vb\otimes\bvb$ as cyclic vector,
\item Filtered tensor product of $\gn[t]$-modules 
$\F_{\zz}\left(
\Pi_{k_1},\cdots,\Pi_{k_p},
\bPi_{\bk_1},\cdots,\bPi_{\bk_{\bp}}\right)$, 
with $\vb\otimes\bvb$ as cyclic vector.  
\end{enumerate}
\end{prop}
\begin{proof}
Set $\bv^{(m)'}=s_{12}\bv^{(m)}$, where
$s_{12}$ denotes the 
reflection with respcect to the simple root $\bep_1-\bep_2$. 
Then we have isomorphisms $\nu,\overline{\nu}$
of $\slt$-modules such that
\be
&
\nu~:~\Pi_{m}\overset{\sim}{\rightarrow}\pi^{(m)}, 
\quad
&\nu\bigl(\exp(e_{23})v^{(m)}\bigr)=u^{(m)},
\\
&\overline{\nu}~:~
\bPi_{m}\overset{\sim}{\rightarrow}\bpi^{(m)}, 
\quad  
&\overline{\nu}\bigl(\exp(e_{23}) \bv^{(m)'}\bigr)
=\bu^{(m)}. 
\en
The filtered tensor product does not change 
if we change the cyclic vector from $\vb\otimes\bvb$ to 
$s_{12}(\vb\otimes\bvb)=\vb\otimes s_{12}(\bvb)$.
Hence the equivalence of (i) and (ii) follows from 
Proposition A.2 in \cite{FJKLM2}.
The equivalence of (ii) and (iii) follows from the PBW 
theorem. 
\end{proof}
If we set $h_{ab}=e_{aa}-e_{bb}$, then in the 
above we have 
\be
&&\nu^{-1}\circ \htil\circ\nu=\frac{1}{3}(h_{12}-h_{23})
+\frac{2}{3}m,
\\
&&\overline{\nu}^{-1}\circ \htil\circ
\overline{\nu}=\frac{1}{3}(h_{12}-h_{23})
-\frac{2}{3}m.
\en

Using the degree operators
$D_0=td/dt$, $D_1=(2h_{12}+h_{23})/3$, 
$D_2=(h_{12}+2h_{23})/3$,
we assign the grading to the corresponding fusion product 
as follows.
\bea
&&\deg e_{12,i}=(i,1,0),
\quad
\deg e_{23,i}=(i,0,1),
\quad
\deg e_{13,i}=(i,1,1),
\label{mixdeg1}\\
&&
\deg \vb\otimes\bvb=
(0,-\frac{|\Mb|+2|\bMb|}{3},-\frac{2|\Mb|+|\bMb|}{3}).
\label{mixdeg2}
\ena
We have then 
\bea
\htil=D_1-D_2+\frac{2}{3}(|\Mb|-|\bMb|).
\label{ht}
\ena

\subsection{Annihilation conditions}\label{subsec:4.3}
The next task is to derive the annihilating conditions for 
the cyclic vector $\vb\otimes \bvb$.
{}For that purpose we introduce 
an abelianization of the fusion product. 

On $U(\gn[t])$ we have a filtration $\{U^{\le i}(\gn[t])\}$
by degrees in $t$. 
Let us consider another filtration $\{U_{\le i}(\gn[t])\}$.
Let $L$ 
be the linear span of $e_{12,i},e_{23,i}$ ($i\ge 0$), 
and set 
\be
&&U_{\le i}(\gn[t])=L U_{\le(i-1)}(\gn[t])
+U_{\le(i-1)}(\gn[t]),
\\
&&U_{\le 0}(\gn[t])=U(\gb[t]),
\quad
U_{\le -1}(\gn[t])=0.
\en
On a cyclic $\gn[t]$-module $W=U(\gn[t])\wb$, we have 
an induced filtration 
\be
&&F^iW=U^{\le i}(\gn[t])\wb,
\\
&&G^iW=U_{\le i}(\gn[t])\wb.
\en
On $\ggr^GW$ the actions of $e_{12,i}$, $e_{23,j}$ are 
commutative. 
Since $\gb[t]\subset \gn[t]$ is central,
$G^iW$ is a $\gb[t]$-module.  

Let $\F$ be the filtered tensor product of $\gn[t]$-modules 
given in Proposition \ref{prop:4.1}, (iii).
The filtration $F$ gives rise to the fusion product
\bea
V_{\Mb,\bMb}={\rm gr}^F\F=
\Pi_{k_1}*\cdots*\Pi_{k_p}*
\bPi_{\bk_1}*\cdots*\bPi_{\bk_{\bp}}. 
\label{fusVV}
\ena
We have 
\bea
\dim V_{\Mb,\bMb}=\sum_{\mb,\bmb}
F_{\Mb,\mb}(1)F_{\bMb,\bmb}(1)\,\chi_{\mb+\bmb}(1,1).
\label{Vmmdim}
\ena

As before, we set 
$e_{ab}(z)=\sum_{i\ge 0}e_{ab,i}z^i$ and 
$\et_{ab}(z)=\prod_{r=1}^{p+\bar{p}}(1-\z_rz)\cdot e_{ab}(z)$. 

\begin{prop}\label{prop:4.2}
The following relations hold on $\ggr^G V_{\Mb,\bMb}$.  
\be
&&\deg_z 
e_{23}(z)^a 
e_{12}(z)^b
e_{13}(z)^c
\vb\otimes\bvb
\\
&&\le
\sum_{i=1}^p\min(a+c,k_i)
+\sum_{i=1}^{\bp}\min(b+c,\bk_i)
-(a+b+c)
\\
&&\qquad \mbox{for any $a, b, c\ge 0$}.
\en
\end{prop}
\begin{proof}
We have $\et_{12}(z)\vb=0$ 
and $\et_{23}(z)\bvb=0$. 
Therefore on $\ggr^G\F$ we obtain 
\bea
&&
\et_{23}(z)^{a} 
\et_{12}(z)^{b} 
\et_{13}(z)^{c} 
\vb\otimes\bvb
\label{et}\\
&&
=\sum_{c_1+c_2=c}
\frac{c!}{c_1!c_2!}
\left(\et_{23}(z)^a \et_{13}(z)^{c_1}\vb\right)
\otimes
\left(\et_{12}(z)^b \et_{13}(z)^{c_2}\bvb\right).
\nn
\ena

At the point $z=\z_1^{-1}$,
$\et_{23}(z)^a\et_{13}(z)^{c_1} \vb$
has a zero of order $\max(c_1+a-k_1,0)$
and 
$\et_{12}(z)^b\et_{13}(z)^{c_2} \bvb$
has a zero of order $c_2+b$. 
Therefore, at $z=\z_1^{-1}$ 
the right hand side of \eqref{et} 
has a zero of order at least $\max(c+a-k_1,0)+b$. 
Proceeding in the same way, we find that 
the left hand side of \eqref{et} is divisible by 
\be
\prod_{i=1}^p(1-\z_iz)^{\max(a+c-k_i,0)+b}
\prod_{i=1}^{\bp}(1-\z_{p+i}z)^{\max(b+c-\bk_i,0)+a}.
\en
Counting degrees and 
passing to $\ggr^F\ggr^G\F=\ggr^G\ggr^F\F=\ggr^GV_{\Mb,\bMb}$, 
we obtain the assertion. 
\end{proof}

\subsection{Subquotient modules and recursion}\label{subsec:4.4}
Let 
\be
\ga=\C X \oplus\C Y \oplus\C Z
\en
be an abelian Lie algebra.
We regard $\ggr^G V_{\Mb,\bMb}$ as an $\ga[t]$-module
where 
$X_i,Y_i,Z_i$ act as $e_{12,i}, e_{23,i},e_{13,i}$,  
respectively. 
Introduce a family of cyclic modules
\bea
W(k_1,\cdots,k_p|\bk_1,\cdots,\bk_{\bp}|l_1,\cdots,l_r)
=U(\ga[t])\bone
\label{mixW}
\ena
by the following defining relations.
\be
&&\deg_z Y(z)^aX(z)^bZ(z)^c \bone
\\
&&
\le 
\sum_{i=1}^p\min(a+c,k_i)
+\sum_{i=1}^{\bp}\min(b+c,\bk_i)
+\sum_{i=1}^r\min(c,l_i)
-(a+b+c)
\en
for all $a,b,c\ge 0$.
 
In the case $p=0$, the module 
$W(|\bk_1,\cdots,\bk_{\bp}|l_1,\cdots,l_r)$
is a special case of the ones studied in \cite{FJKLM2}. 

{}From Proposition \ref{prop:4.2}, 
we have also a surjection of $U(\ga[t])$-modules 
\bea
W(k_1,\cdots,k_p|\bk_1,\cdots,\bk_{\bp}|)
\longrightarrow
{\rm gr}^G V_{\Mb,\bMb}
\longrightarrow 0. 
\label{surj2}
\ena

The rest of the working is entirely similar to the previous 
section. 
Set
\be
&&
W=W(k_1,\cdots,k_p|\bk_1,\cdots,\bk_{\bp}|l_1,\cdots,l_r),
\\
&&
W'=W(k_2,\cdots,k_p|\bk_1,\cdots,\bk_{\bp}|k_1,l_1,\cdots,l_r),
\\
&&
W''=W(k_1-1,k_2,\cdots,k_p|\bk_1,\cdots,\bk_{\bp}|l_1,\cdots,l_r).
\en
Let further $\bW'$ be the subspace 
of $W$ generated by $\bone$ over 
$X_i,Z_i$ ($i\ge 0$) and $Y_i$ ($i\ge 1$). 

\begin{thm}\label{thm:4.3}
Suppose $k_1\ge k_2,\cdots,k_p$. 
Then there exists surjective maps of 
$\C[Z_i~(i\ge 0)]$-modules 
\be
&&\iota~:~W'\rightarrow \bW',
\\
&&\phi~:~W''\rightarrow W/\bW'.
\en
The maps are
\be
&&
\iota(X_i)=X_i,\quad \iota(Y_i)=Y_{i+1},\quad \iota(Z_i)=Z_i,
\\
&&
\phi(b Y_0^{(m)})=bY_0^{(m+1)},
\en
where $b$ is not divisible by $Y_0$.
Similar maps exist if we exchange the roles of 
$k_1,\cdots,k_p$ with $\bk_1,\cdots,\bk_{\bp}$
and $Y_i$ with $X_i$.
\end{thm}

We have a $\Z_{\ge 0}^4$ grading on 
\eqref{mixW} given by 
\be
&&\deg X_i=(i,1,0,0),
\quad
\deg Y_i=(i,0,1,0),
\quad
\deg Z_i=(i,0,0,1),
\\
&&\deg {\bf 1}=(0,0,0,0).
\en
The character satisfies the recursive estimate
\be
&&\ch_{q,z_1,z_2,z_3}
W(k_1,\cdots,k_p|\bk_1,\cdots,\bk_{\bp}|l_1,\cdots,l_r)
\\
&&
\le\ch_{q,z_1,qz_2,z_3}
W(k_2,\cdots,k_p|\bk_1,\cdots,\bk_{\bp}|k_1,l_1,\cdots,l_r)
\\
&&+
z_2\ch_{q,z_1,z_2,z_3}
W(k_1-1,k_2,\cdots,k_p|\bk_1,\cdots,\bk_{\bp}|l_1,\cdots,l_r),
\en
under the assumption that $k_1\ge k_2,\cdots,k_p$. 
{}From this and the known initial condition for $p=0$, 
we obtain 
\be
&&\ch_{q,z_1,z_2,z_3}
W(k_1,\cdots,k_p|\bk_1,\cdots,\bk_{\bp}|l_1,\cdots,l_r)
\\
&&\le\sum_{\mb,\bmb}
F_{\Mb,\mb}(q)F_{\bMb,\bmb}(q)
z_2^{|\Mb|-|\mb|}z_1^{|\bMb|-|\bmb|}\chi_{\mb+\bmb+\nb}(q,z_3),
\en
where $\nb=(n_1,\cdots,n_k)$, 
$n_a=\sharp\{j\mid l_j=a\}$.
{}From \eqref{surj2} and \eqref{Vmmdim}, 
we have an equality for $r=0$. 

In order to obtain the character of the fusion product, 
we specialize $z_3=z_1z_2$ and supply an overall 
power $z_1^{-(|\Mb|+2|\bMb|)/3}z_2^{-(2|\Mb|+|\bMb|)/3}$, 
to take into account the degree of the cyclic vector
\eqref{mixdeg2}. 

\begin{thm}\label{thm:sl3}
There exists a filtration $H$ of $\gr^GV_{\Mb,\bMb}$ 
by $\C[Z_i~(i\ge 0)]$-modules such that
\be
\ggr^H\ggr^G V_{\Mb,\bMb}=\bigoplus_{\mb,\bmb}
\mathcal{M}_{\Mb\mb;\bMb\bmb}\otimes\pi_{\mb+\bmb},
\en
where $\mathcal{M}_{\Mb\mb;\bMb\bmb}$ is a trivial module with
the character
\be
\ch_q\mathcal{M}_{\Mb\mb;\bMb\bmb}=
F_{\Mb,\mb}(q)\,F_{\bMb,\bmb}(q). 
\en
The following formula holds for the character of
the fusion product 
\be
&&
\ch_{q,z_1,z_2}V_{\Mb,\bMb}
=\sum_{\mb, \bmb}
F_{\Mb,\mb}(q)\,F_{\bMb,\bmb}(q)
\,\ch_{q,(z_1z_2)}\pi_{\mb+\bmb}
\\
&&\quad\times
(z_1^{-1}z_2)^{(|\Mb|-|\bMb|)/3-(|\mb|-|\bmb|)/2}.
\en
\end{thm}

{}Finally Theorem \ref{thm:2.2} follows by taking the quotient
with respect to $Z_0$, $Z_1^{k-l+1}$,
$D_2+l/2$ and using the 
information about the dimension \eqref{dimmix}.
In view of \eqref{mixdeg2} and \eqref{ht}, 
the character is obtained by settin 
$z_1=z,z_2=z^{-1}$ and multiplying by 
$z^{2(|\Mb|-|\bMb|)/3}$. 
\bigskip

\noindent
{\it Acknowledgments.}\quad
BF is partially supported by the grants 
RFBR-02-01-01015, RFBR-01-01-00906 and INTAS-00-00055. 
SL is partially supported by the grants RFBR-02-01-01015 and 
RFBR-01-01-00546.
JM is partially supported by 
the Grant-in-Aid for Scientific Research (B2) no.12440039, 
and TM is partially supported by 
(A1) no.13304010, Japan Society for the Promotion of Science.


\begin{thebibliography}{[FKLMM3]}

\bibitem[FF]{FF}
B.~L. Feigin and E.~Feigin,
\newblock $q$-characters of the tensor products in
$\slt$-case,
\newblock math.QA/0201111 (2002).

\bibitem[FJKLM1]{FJKLM1}
B.~Feigin, M.~Jimbo, R.~Kedem, S.~Loktev and T.~Miwa,
Spaces of coinvariants and fusion product I. 
{}From equivalence theorem to Kostka polynomials.
\newblock math.QA/0205324

\bibitem[FJKLM2]{FJKLM2}
B.~Feigin, M.~Jimbo, R.~Kedem, S.~Loktev and T.~Miwa,
Spaces of coinvariants and fusion Product II. 
{}$\slth$ character formulas in terms of Kostka polynomials.
\newblock math.QA/0208156

\bibitem[FKLMM1]{FKLMM1}
B.~Feigin, R.~Kedem, S.~Loktev, T.~Miwa and E.~Mukhin,
\newblock Combinatorics of the $\slth$ Spaces of Coinvariants,
\newblock Transformation Groups {\bf 6} (2001) 25--52.

\bibitem[FKLMM2]{FKLMM2}
B.~Feigin, R.~Kedem, S.~Loktev, T.~Miwa and E.~Mukhin,
\newblock Combinatorics of the $\slth$ Spaces of Coinvariants:
Loop Heisenberg modules and recursion, 
\newblock math.QA/0009198


\bibitem[FKLMM3]{FKLMM3}
B.~Feigin, R.~Kedem, S.~Loktev, T.~Miwa and E.~Mukhin,
\newblock Combinatorics of the $\slth$ spaces of coinvariants:
Dual functional realization and recursion,
 \newblock math.QA/0012190

\bibitem[FL]{FL}
B.~L. Feigin and S.~Loktev,
\newblock On generalized Kostka polynomials and
quantum Verlinde rule,
math.QA/9812093,
\newblock {\em Amer.~Math.~Sci.~Transl.}
{\bf 194} (1999) 61--79.

\bibitem[SS]{SS}
A.~Schilling and M.~Shimozono, 
\newblock Fermionic formulas for level-restricted generalized
Kostka polynomials and coset branching functions, 
math.QA/0001114, 
\newblock {\em Commun. Math. Phys.} {\bf 220} (2001) 105--164. 

\bibitem[SW]{SW}
A.~Schilling and S.~O.~Warnaar,
\newblock Inhomogeneous lattice paths, 
generalized Kostka polynomials 
and $A_{n-1}$ supernomials, 
\newblock math.QA/9802111, 
\newblock {\em Commun. Math. Phys.} {\bf 202} (1999) 359--401.

\end{thebibliography}
\end{document}